\documentclass{aims}
\usepackage{amsmath}
  \usepackage{paralist}
  \usepackage{graphics} 
  \usepackage{epsfig} 
\usepackage{graphicx}  \usepackage{epstopdf}
 \usepackage[colorlinks=true]{hyperref}
\hypersetup{urlcolor=blue, citecolor=red}

\theoremstyle{definition}

\title[Note on the ideal exponentiation] 
      {Note on the ideal exponentiation in imaginary quadratic number and function fields}

\author[Soufiane Mezroui]{}


 \email{mezroui.soufiane@yahoo.fr}

\begin{document}

\begin{abstract}
By using a new formula of cubing ideals in imaginary quadratic number and function fields combined with Shank's NUCOMP algorithm, Imbert et al. presented a fast algorithms that compute a reduced output of cubing ideals and keep the sizes of the intermediate operands small. The authors asked whether there are
specialized formulas to compute higher powers of ideals. In this note, we will derive such formulas. 
\end{abstract}

\maketitle

\centerline{\scshape Soufiane Mezroui}
\medskip
{\footnotesize
 \centerline{Abdelmalek Essadi University}
   \centerline{National School of Applied Sciences of Tangier (ENSAT)}
   \centerline{Laboratory of Information and Communication Technologies (LabTIC)}
   \centerline{ENSA Tanger, Route Ziaten, BP 1818,}
   \centerline{Tanger principale, Morocco}
} 



\bigskip


\section{Introduction}

Let $\mathbb{K}=\mathbb{F}_q(C)=\mathbb{F}_q(x,y)$ be an imaginary quadratic function field of genus $g$ over a finite field $\mathbb{F}_q$, where $C$ is an imaginary hyperelliptic curve defined by
\begin{equation}
C: y^{2}+h(x)y=f(x),
\end{equation}
such that $f,h\in \mathbb{F}_q[x]$, $f$ is monic of degree $2g+1$, and $h=0$ if $q$ is odd and is monic with $deg(h)\leq g$ if $q$ is even. 

Remember that the maximal order $\mathbb{F}_q[x,y]$ of $\mathbb{F}_q(x, y)$ is an integral domain, it is also a $\mathbb{F}_q[x]$-module of rank $2$ with $\mathbb{F}_q[x]$-basis $\{1,y\}$. The non-zero integral ideals in $\mathbb{F}_q[x, y]$ are the $\mathbb{F}_q[x]$-modules of the following form $\mathfrak{a}= \mathbb{F}_q[x] S Q + \mathbb{F}_q[x] S(P + y)$ where $P, Q, S\in \mathbb{F}_q[x]$ and $Q$ divides $f + hP - P^2$. Here, $S$ and $Q$ are unique up to factors in $\mathbb{F}_{q}*$ and $P$ is unique modulo $Q$. For brevity, we will write $\mathfrak{a}=S(Q, P)$. 

By using a new formula of cubing ideals in imaginary quadratic number and function fields combined with Shank's NUCOMP algorithm \cite{Jacobson,Shanks}, Imbert et al. \cite{imbert} presented a fast algorithms that compute a reduced output of $\mathfrak{a}^{3}$ and keep the sizes of the intermediate operands small. The authors asked whether there are specialized formulas to compute higher powers of ideals. In this note, we will derive such formulas.

\section{Deriving a formula to compute higher powers of ideals in imaginary quadratic function fields}

Suppose that $S=1$ and $gcd(Q,2P-h)=1$ satisfying $u_{1}Q+v_{1}(2P-h)=1$, $u_{1},v_{1}\in\mathbb{F}_q[x]$. We will show that 
\begin{equation}
\mathfrak{a}^{n+1}=(Q^{n+1},P+v_{1}Q R S_{n+1}), 
\end{equation}
where $R=\frac{f+Ph-P^2}{Q}$ and $S_{n+1}$ is computed recursively by 
$$
S_{n+1}=1+(1+v_{1}h-2P v_{1})S_{n}-v_{1}^{2}Q R S_{n}^{2},
$$
such that $S_{1}=0$. 

Indeed, assume that $\mathfrak{a}^{n}=(Q^{n},P+v_{1}Q R S_{n})$. Using the formula to compute the product of ideals (see \cite[pp.243]{imbert}) we obtain 
\begin{equation}
\mathfrak{a}^{n+1}=\mathfrak{a}\times \mathfrak{a}^{n}=(Q,P)\times (Q^{n},P+v_{1}Q R S_{n}),
\end{equation}
such that 
\begin{equation} 
\begin{split}
gcd(Q,Q^{n},2P-h+v_{1}Q R S_{n})&=gcd(Q,2P-h+v_{1}Q R S_{n})\\
&=gcd(Q,2P-h)\\
&=1.
\end{split}
\end{equation} 
 Let be $U,V,W\in\mathbb{F}_q[x]$ satisfying
\begin{equation} \label{eq1}
\begin{split}
1 & = U Q+V Q^{n}+W (2P-h+v_{1}Q R S_{n}) \\
 & = (U+W v_{1} R S_{n}) Q+V Q^{n}+W (2P-h), 
\end{split}
\end{equation}
then $U+W v_{1} R S_{n}=u_1$, $V=0$ and $W=v_1$. Put $P_{n}=P+v_{1}Q R S_{n}$, we deduce the following
\begin{equation} 
\begin{split}
\mathfrak{a}^{n+1} & = \mathfrak{a}\times \mathfrak{a}^{n} \\
 & = (Q,P)\times (Q^{n},P+v_{1}Q R S_{n}), \\
 & = \left(Q^{n+1},P+v_{1}Q R S_{n}+Q^{n}\left(V(-v_{1}Q R S_{n})+W \frac{f+P_{n}h-P_{n}^2}{Q^{n}}\right)\right),\\
 & = \left(Q^{n+1},P+v_{1}Q R S_{n}+Q^{n} v_{1} \frac{f+P_{n}h-P_{n}^2}{Q^{n}}\right),\\
 & = (Q^{n+1},P+v_{1}Q R S_{n}+ v_{1} (f+(P+v_{1}Q R S_{n})h-(P+v_{1}Q R S_{n})^2)),\\
 & = (Q^{n+1},P+v_{1}Q R S_{n}\\
&\,\,\,+v_{1} (f+P h-P^{2}+v_{1}Q R S_{n}h-2 P v_{1}Q R S_{n}-v_{1}^{2}Q^{2} R^{2} S_{n}^{2})),\\
 & = (Q^{n+1},P+v_{1}Q R (1+(1+v_{1}h-2P v_{1})S_{n}-v_{1}^{2}Q R S_{n}^{2})).
\end{split}
\end{equation}
Then $\mathfrak{a}^{n+1}=(Q^{n+1},P+v_{1}Q R S_{n+1})$ with 
$$
S_{n+1}=1+(1+v_{1}h-2P v_{1})S_{n}-v_{1}^{2}Q R S_{n}^{2},
$$ 
which gives the result. Notice that $S_{1}=0$, $S_{2}=1$, and $S_{3}=2+v_{1}(h-2P-v_{1}Q R)$ which is the formula obtained in \cite{imbert}. 

\section{The imaginary quadratic number fields case}
Let $\mathbb{Q}(\sqrt{\Delta})$ be an imaginary quadratic number fields and let 
$$
\mathfrak{a}= \mathbb{Z}S Q + \mathbb{Z} S\frac{P + \sqrt{\Delta}}{2},
$$
be an integral ideal such that $P,Q,S\in\mathbb{Z}$ and $4Q\mid P^{2}-\Delta$. We will denote it for brievty by $\mathfrak{a}=S(P,Q)$. 

Similarly as in the function fields case, suppose $S=1$ and $gcd(Q,P)=1$ satisfying $u_{1}Q+v_{1} P=1$ where $u_{1},v_{1}\in\mathbb{Z}$. The formula to compute higher powers of ideals $\mathfrak{a}^{n+1}$ is given by
$$
\mathfrak{a}^{n+1}=(Q^{n+1},P-2 Q S_{n+1}),
$$
where $S_{n+1}$ is obtained recursively by 
$$
S_{n+1}=c v_{1}+Q u_{1} S_{n}+Q v_{1} S_{n}^{2},
$$
such that $S_{1}=0$ and $c=\frac{P^{2}-\Delta}{4 Q}$.




\end{document}